\documentclass[12pt,reqno]{amsart}
\usepackage{amssymb, amsmath}
\newtheorem{thm}{Theorem}

\numberwithin{defn}{section}
\numberwithin{thm}{section}
\numberwithin{Lemma}{section}
\numberwithin{Corollary}{section}
\numberwithin{Example}{section}
\numberwithin{subsection}{section}
\numberwithin{Remark}{section}
\numberwithin{equation}{section}
\numberwithin{ppn}{section}
\begin{document}
\title[ New Efficient Steffensen Type method  for Solving. . .]
{ New Efficient Steffensen Type method  for Solving Nonlinear Equations}
\author{J. P. Jaiswal}
\date{}
\maketitle


\textbf{Abstract.} 
In the present paper, by approximating the derivatives in the Kou et al. \cite{Kou} fourth-order method by central difference quotient, we obtain new modification of this method free from derivatives. We prove the important fact that the method obtained preserve their order of convergence, without calculating any derivative. Finally, numerical tests confirm that our method give the better performance as compare to other well known Steffensen type methods.
\\



\section{Introduction}
A large number of papers have been written about iterative methods for the solution of the nonlinear equations. In this paper, we consider the problem of finding a simple root $x^*$ of a function  $f: D\subseteq \Re \longrightarrow \Re$, i.e. $f(x^*)=0$ and $f'(x^*)\neq0$. The famous Newton's method for finding $x^*$ uses the iterative method 
\begin{equation*}
x_{n+1}=x_{n}-\frac{f(x_n)}{f'(x_n)},
\end{equation*}
starting from some initial value $x_0$. The Newton's method is an important and basic method which converges quadratically in some neighborhood of simple root $x^*$. However, when the first order derivative of the function $f(x)$ is unavailable or is expensive to compute, the Newton's method is still restricted in practical applications. In order to avoid computing the first order derivative, Steffensen in \cite{Kincaid} proposed the following derivative-free method. It is well known that the forward-difference approximation for $f'(x)$ at $x$ is 
\begin{equation*}
f'(x)\approx \frac{f(x+h)-f(x)}{h}.
\end{equation*}
If the derivative $f'(x_n)$ is replaced by the forward-difference approximation with $h=f'(x_n)$ i.e.
\begin{equation*}
f'(x_n)\approx \frac{f(x_n+f(x_n))-f(x)}{f(x_n)},
\end{equation*}
the Newton's method becomes
\begin{equation*}
x_{n+1}=x_{n}-\frac{f(x_n)^2}{f(x_n+f(x_n))-f(x_n)},
\end{equation*}
which is the famous Steffensen's method \cite{Kincaid}. The Steffensen's method is based on forward-difference approximation to derivative. This method is a tough competitor of Newton's method. Both the methods are quadratic  convergence, both require two functions evaluation per iteration but Steffensen's method is derivative free. The idea of removing derivatives from the iteration process is very significant. Recently, many high order derivative-free methods are built according to the Steffensen's method, see [\cite{Cordero}, \cite{Liu}, \cite{Ren}, \cite{Zheng}, \cite{Petkovic}, \cite{Jain}, \cite{Dehghan1}, \cite{Dehghan2}] and the references therein.\

Jain \cite{Jain} presented the third-order method (Jain Method): 
\begin{eqnarray}\label{eqn:11}
y_n&=&x_n-\frac{f(x_n)^2}{f(x_n+f(x_n))-f(x_n)},\nonumber \\
x_{n+1}&=&x_n-\frac{f(x_n)^3}{\{f(x_n+f(x_n))-f(x_n)\}\{f(x_n)-f(y_n)\}}.
\end{eqnarray}

Dehghan et al. \cite{Dehghan1} proposed two third-order Steffensen type method (Dehghan Method I):
\begin{eqnarray}\label{eqn:12}
y_n&=&x_n-\frac{2f(x_n)^2}{f(x_n+f(x_n))-f(x_n-f(x_n))},\nonumber \\
x_{n+1}&=&x_n-\frac{2f(x_n)\{f(x_n)+f(y_n)\}}{f(x_n+f(x_n))-f(x_n-f(x_n))}.
\end{eqnarray}
and (Dehghan Method II)
\begin{eqnarray}\label{eqn:13}
y_n&=&x_n+\frac{2f(x_n)^2}{f(x_n+f(x_n))-f(x_n-f(x_n))},\nonumber \\
x_{n+1}&=&x_n-\frac{2f(x_n)\{f(y_n)-f(x_n)\}}{f(x_n+f(x_n))-f(x_n-f(x_n))}.
\end{eqnarray}

Again Dehghan et al. \cite{Dehghan2} introduced a new third-order Steffensen type method (Dehghan Method III):
\begin{eqnarray}\label{eqn:14}
y_n&=&x_n+\frac{2f(x_n)^2}{f(x_n+f(x_n))-f(x_n-f(x_n))},\nonumber \\
x_{n+1}&=& x_n-\frac{2f(x_n)}{f(y_n)f(u)+f(x_n)f(v)},
\end{eqnarray}
where $f(u)=f(x_n+f(x_n))-f(x_n-f(x_n))$ \\
and $f(v)=f(y_n+f(y_n))-f(y_n-f(y_n))$.\\

Recently Cordero et al. \cite{Cordero} presented a fourth-order Steffensen type method (Cordero Method):
\begin{eqnarray}\label{eqn:15}
y_n&=&x_n-\frac{2f(x_n)^2}{f(x_n+f(x_n))-f(x_n-f(x_n))},\nonumber \\
x_{n+1}&=&x_n-\frac{2f(x_n)^2}{f(x_n+f(x_n))-f(x_n-f(x_n))}.\frac{f(y_n)-f(x_n)}{2f(y_n)-f(x_n)}.
\end{eqnarray}
Other  Steffensen type methods and their applications are discussed in [\cite{Liu}, \cite{Ren}, \cite{Zheng}, \cite{Petkovic}]. 

The purpose of this paper is to develop a new fourth-order derivative-free method and give the convergence analysis. This paper is organized as follows. In Section 2, we present a new two-step fourth-order iterative method for solving nonlinear equations. In this method we approximate the derivative of the function by central difference quotient. The new method is free from derivative. We prove that the order of convergence of the new method is four. Numerical examples show better performance of our method in section 4. Section 5 is a short conclusion.




\section{Development of the method and analysis of convergence}
Let us consider the fourth-order of method proposed by Kou et al. in \cite{Kou}:
\begin{eqnarray}\label{eqn:12a}
y_n&=&x_n-\frac{f(x_n)}{f'(x_n)},\nonumber \\
x_{n+1}&=&x_n-\theta\frac{f(x_n)+f(y_n)}{f'(x_n)}-(1-\theta)\frac{f(x_n)^2}{f'(x_n)[f(x_n)-f(y_n)]},
\end{eqnarray}
where $\theta\in R$. The convergence order of these methods is at least three, and for  $\theta=-1$ the order is four. Specifically, the third order method obtained taking $\theta=0$ in $(\ref{eqn:12a})$, and the forth order method corresponding to $\theta=-1$. Here we are going to use this fourth-order method of this family.
We try to use forward-difference approximation of the derivative appeared in $(\ref{eqn:12a})$, then the the order of convergence of the new method goes down to three and error expression for it is $(-c_2^2-c_2^2/c_1)e_n^3+O(e_n^4)$, where $c_i$ is defined later. For this reason, we have used the central-difference approximation in $(\ref{eqn:12a})$, obtaining a modified Kou et al. method that preserves the convergence order and is derivative free. 

By using a symmetric difference quotient
\begin{equation*}
f'(x_n)\approx \frac{f(x_n+f(x_n))-f(x_n-f(x_n))}{2f(x_n)},
\end{equation*}

By approximating the derivative by the central-difference in $(\ref{eqn:12a})$ for $\theta=-1$, we obtain a new method free from derivatives, that we call the \textsl{modified Kou method free from derivative (MKDF)}:
\begin{eqnarray}\label{eqn:22}
y_n&=&x_n-\frac{2f(x_n)^2}{f(x_n+f(x_n))-f(x_n-f(x_n))},\nonumber \\
x_{n+1}&=&x_n+\frac{2f(x_n)\{f(x_n)+f(y_n)\}}{f(x_n+f(x_n))-f(x_n-f(x_n))}\nonumber \\
        &&-\frac{4f(x_n)^3}{\{f(x_n+f(x_n))-f(x_n-f(x_n))\}\{f(x_n)-f(y_n)\}},
\end{eqnarray}
Now we are going to prove the method \textsl{MKDF} have orders of convergence four.
\begin{thm}
Let $x^* \in I$ be a simple zero of a sufficiently differentiable function $f:I\subseteq \Re \rightarrow \Re$ in an open interval I. If $x_0$ is sufficiently close to $x^*$, then the modified Kou method free from derivative defined by $(\ref{eqn:22})$ has order of convergence four and satisfies the error equation $(\ref{eqn:29})$.
\end{thm}
\begin{proof}
By applying the Taylor series expansion theorem and taking account $f(x^*)=0$, we can write 
\begin{equation}\label{eqn:23}
f(x_n)=c_1e_n+c_2e_n^2+c_3e_n^3+c_4e_n^4+c_5e_n^5+c_6e_n^6+c_7e_n^7+c_8e_n^8+O(e_n^9),
\end{equation}
where $c_k=\frac{f^{k}(x^*)}{\lfloor k},\ k=1, 2, . . .$ and $e_n$ be the error in $x_n$ after $n$ iterations  i.e. $e_n=x_n-x^*$;
\begin{eqnarray}\label{eqn:24}
f(x_n+f(x_n))&=&(c_1^2 + c_1)e_n+(c_2c_1^2 + 3c_2c_1 + c_2)e_n^2 \nonumber\\
            &+&(c_3c_1^3 + 3c_3c_1^2 + 2c_1c_2^2 + 4c_3c_1 + 2c_2^2 + c_3)e_n^3 \nonumber\\   
            &+&(c4 + c_2(c_2^2 + 2c_1c_3) + 5c_1c_4 + 5c_2c_3 + 6c_1^2c_4 \nonumber\\
            &&+ 4c_1^3c_4 + c_1^4c_4 + 6c_1c_2c_3+3c_1^2c_2c_3)e_n^4 \nonumber\\
            &+&O(e_n^5).
\end{eqnarray}
and 
\begin{eqnarray}\label{eqn:25}
f(x_n-f(x_n))&=&(-c_1^2 + c_1)e_n+(c_2c_1^2-3c_2c_1 + c_2)e_n^2 \nonumber\\
            &+&(-c_3c_1^3 + 3c_3c_1^2 + 2c_1c_2^2 - 4c_3c_1 - 2c_2^2 + c_3)e_n^3 \nonumber\\   
            &+&(c4 + c_2(c_2^2 + 2c_1c_3) -5c_1c_4 - 5c_2c_3 + 6c_1^2c_4 \nonumber\\
            &&- 4c_1^3c_4 + c_1^4c_4 + 6c_1c_2c_3-3c_1^2c_2c_3)e_n^4 \nonumber\\
            &+&O(e_n^5).
\end{eqnarray}
Further more it can be easily find 
\begin{eqnarray}\label{eqn:26}
&&\frac{2.f(x_n)^2}{f(x_n+f(x_n))-f(x_n-f(x_n))}\nonumber\\
=&&e_n+(-c_2/c_1)e_n^2 + ((2c_2^2)/c_1^2 -(2c_3)/c_1- c_1c_3)e_n^3\nonumber\\     
 &&+(c_2c_3-4c_1c_4-(3c_4)/c_1 - (4c_2^3)/c_1^3+(7c_2c_3)/c_1^2)e_n^4 \nonumber\\
 &&+O(e_n^5) .
\end{eqnarray}
By considering this relation and expression of $y_n$ in the equation $(\ref{eqn:22})$, we obtain
\begin{eqnarray}\label{eqn:27}
y_n&=&x^*+(c_2/c_1)e_n^2-((2c_2^2)/c_1^2 -(2c_3)/c_1- c_1c_3)e_n^3\nonumber\\     
   &&-(c_2c_3-4c_1c_4-(3c_4)/c_1 - (4c_2^3)/c_1^3+(7c_2c_3)/c_1^2)e_n^4 \nonumber\\
  &&+O(e_n^5) .
\end{eqnarray}
At this time, we should expand $f(y_n)$ around the root by taking into consideration $(\ref{eqn:27})$. Accordingly, we have
\begin{eqnarray}\label{eqn:28}
f(y_n)&=&c_2e_n^2+c_1(c_1c_3+(2c_3)/c_1 - (2c_2^2)/c_1^2)e_n^3 \nonumber\\     
   &&+(c_1(4c_1c_4 - c_2c_3 + (3c_4)/c_1 + (4c_2^3)/c_1^3 - (7c_2c_3)/c_1^2) + c_2^3/c_1^2)e_n^4\nonumber\\
  &&+O(e_n^5) .
\end{eqnarray}  
By using $(\ref{eqn:23})$, $(\ref{eqn:28})$, $(\ref{eqn:24})$ and $(\ref{eqn:25})$ in the last expression of $(\ref{eqn:22})$, we obtain
\begin{eqnarray}\label{eqn:29}
e_{n+1}&=&(-4c_1c_4-3c_2c_3-(11c_4)/c_1-(2c_2^3)/c_1^3-(17c_2c_3)/c_1^2)e_n^4\nonumber\\
        &&+O(e_n^5). 
\end{eqnarray}

\end{proof}

\newpage
\section{Numerical Tests}
In this section, in order to compare the our new method with Steffensen method, Jain method, Dehghan method I, Dehghan method II, Dehghan method III and Cordero method, we give some numerical examples. For this consider the following functions: 
 
\begin{table}[htb]
 \caption{ Test functions and their roots.}
  \begin{tabular}{ll} \hline
Non-linear functions   			     & \hspace{50pt}Roots \\  \hline
$f_1(x)=\sin^2(x)-x^2+1$   			& \hspace{50pt}1.404492 \\ 
$f_2(x)=x^2-e^x-3x+2$      			& \hspace{50pt}0.257530 \\ 
$f_3(x)=cos(x)-x$          			& \hspace{50pt}0.739085 \\ 
$f_4(x)=\cos(x)-xe^x+x^2$  			& \hspace{50pt}0.639154 \\ 
$f_5(x)=e^x-1.5-\arctan(x)$   	& \hspace{50pt}0.767653 \\ 
$f_6(x)=8x-\cos(x)-2x^2$ 				& \hspace{50pt}0.128077 \\ 
$f_7(x)=\ln(x^2+x+2)-x+1$ 			& \hspace{50pt}4.152590 \\ \hline 
  \end{tabular}
  \label{tab:abbr}
\end{table}
Table 2-8 shows the comparison of these methods for these functions. All the numerical computations have been carried out using MATHEMATICA 8. The numerical results show that the our proposed method is efficient.
 \begin{table}[htb]
 \caption{Errors Occurring in the estimates of the root of function $f_1(x)=\sin^2(x)-x^2+1$ after third iteration by the method described with initial guess $x_0=1$.}
  \begin{tabular}{llll} \hline
Methods &  $\left|f_1(x_3)\right|$  \\ \hline
$Steffensen \ Method$   &0.27307e-3\\ 
$Jain \ Method$        &0.15974e-11 \\
$Dehghan \ I\  Method$    &0.12208e-3 \\
$Dehghan \ II \ Method$   &0.27058e-7 \\  
$Dehghan \ III \ Method$  &0.23473e-9 \\
$Cordero \ Method$      &0.16292e-9 \\
$MKDF \ Method$         &0.47200e-25 \\ 
\hline   
  \end{tabular}
  \label{tab:abbr}
\end{table}                                                 
 
 \begin{table}[htb]
 \caption{Errors Occurring in the estimates of the root of function $f_2(x)=x^2-e^x-3x+2$ by the method described with initial guess $x_0=0.7$.}
  \begin{tabular}{llll} \hline
Methods & $\left|f_2(x_3)\right|$ \\\hline
$Steffensen \ Method$   &0.79151e-6\\
$Jain \ Method$        &0.99771e-31 \\
$Dehghan \ I \ Method$    &0.27162e-5 \\
$Dehghan \ II \ Method$   &0.20573e-22 \\
$Dehghan \ III \ Method$  &0.70112e-12 \\
$Cordero$      &0.89687e-6 \\
$MKDF \ Method$         &0.77899e-31 \\
 \hline  
  \end{tabular}
  \label{tab:abbr}
\end{table} 

 \newpage 
 \begin{table}[htb]
 \caption{Errors Occurring in the estimates of the root of function $f_3(x)=cos(x)-x$ by the method described with initial guess $x_0=1$.}
  \begin{tabular}{llll} \hline
Methods & $\left|f_3(x_3)\right|$ \\\hline
$Steffensen \ Method$   &0.82149e-10\\
$Jain \ Method$        &0.10235e-34 \\
$Dehghan\ I \ Method$    &0.17099e-23 \\
$Dehghan\ II\ Method$   &0.61489e-33 \\ 
$Dehghan \ III \ Method$  &0.19366e-28 \\
$Cordero \ Method$      &0.14803e-15 \\
$MKDF \ Method$         &0.10720e-64 \\
   \hline
  \end{tabular}
  \label{tab:abbr}
\end{table} 

 
 \begin{table}[htb]
 \caption{Errors Occurring in the estimates of the root of function $f_4(x)=\cos(x)-xe^x+x^2$ by the method described with initial guess $x_0=1$.}
  \begin{tabular}{llll} \hline
Methods & $\left|f_4(x_3)\right|$ \\\hline
$Steffensen \ Method$   &0.46704e-2\\
$Jain \ Method$        &0.15915e-10 \\
$Dehghan \ I \ Method$    &0.85626e-2 \\ 
$Dehghan \ II\  Method$   &0.40435e-7 \\ 
$Dehghan \ III \ Method$  &0.15970e-5 \\
$Cordero \ Method$      &0.19461e-2 \\
$MKDF \ Method$         &0.16214e-27 \\
   \hline
  \end{tabular}
  \label{tab:abbr}
\end{table} 
 
 
 \begin{table}[htb]
 \caption{Errors Occurring in the estimates of the root of function $f_5(x)=e^x-1.5-\arctan(x)$ by the method described with initial guess $x_0=1$.}
  \begin{tabular}{llll} \hline
Methods & $\left|f_5(x_3)\right|$ \\\hline
$Steffensen \ Method$   &0.78069e-3\\
$Jain \ Method$        &0.99920e-15 \\
$Dehghan \ I \ Method$    &0.43288e+0 \\
$Dehghan \ II \ Method$   &0.11102e-15 \\
$Dehghan \ III \ Method$  &0.21494e-10 \\
$Cordero Method$      &0.55888e-13 \\
$MKDF \ Method$         &0.11102e-15 \\
   \hline
  \end{tabular}
  \label{tab:abbr}
\end{table} 

\newpage  
 \begin{table}[htb]
 \caption{Errors Occurring in the estimates of the root of function $f_6(x)=8x-\cos(x)-2x^2$ by the method described with initial guess $x_0=1$.}
  \begin{tabular}{llll} \hline
Methods & $\left|f_6(x_3)\right|$ \\\hline
$Steffensen \ Method$   &0.77299e+1\\
$Jain \ Method$        &0.37235e-3 \\ 
$Dehghan \ I \ Method$    &0.56299e-9 \\ 
$Dehghan \ II \ Method$   &0.29397e-2 \\
$Dehghan \ III \ Method$  &0.52774e-14 \\
$Cordero \ Method$      &0.77803e-24 \\
$MKDF \ Method$         &0.38994e-14 \\
  \hline
  \end{tabular}
  \label{tab:abbr}
\end{table} 

 
 \begin{table}[htb]
 \caption{Errors Occurring in the estimates of the root of function $f_7(x)=\ln(x^2+x+2)-x+1$ by the method described with initial guess $x_0=3.6$.}
  \begin{tabular}{llll} \hline
Methods & $\left|f_7(x_3)\right|$ \\\hline
$Steffensen \ Method$   &0.51451e-13\\
$Jain \ Method $        &0.11061e-42 \\ 
$Dehghan \ I \ Method$    &0.15144e-29 \\ 
$Dehghan \ II \ Method$   &0.48976e-33 \\
$Dehghan \ III \ Method$  &0.19986e-49 \\
$Cordero \ Method$      &0.43626e-79 \\
$MKDF \ Method$         &0.39633e-74 \\ \hline
  \end{tabular}
  \label{tab:abbr}
\end{table} 
                              
\section{Conclusion} 
The problem of locating roots of nonlinear equations occurs frequently in scientific and engineering works. In this paper, we have introduced a new derivative free fourth-order method for solving nonlinear equations.  By theoretical result, we confirm that the new method archive fourth order convergence. This method can be used for solving nonlinear equations without computing derivatives. The new method is also compared in the performance with some well known steffensen type methods. Numerical results indicate that the our proposed method gives better performance.


\textsc{Jai Prakash Jaiswal\\
Department of Mathematics, \\
Maulana Azad National Institute of Technology,\\
Bhopal, M.P., India-462051}.\\
E-mail: { asstprofjpmanit@gmail.com; jaiprakashjaiswal@manit.ac.in}.
 

\end{document}